\documentclass{article}
\usepackage{amsmath,amssymb,bbm,theorem,stmaryrd, multicol}
\usepackage[french]{babel}
\usepackage[utf8]{inputenc}
\usepackage{hyperref}

\input xy
\xyoption{all}

\def\KK {\mathbbm K}

\def\QQ {\mathbbm Q}

\def\CC {\mathbbm C}

\def\Ker {\mathrm{Ker}}

\newtheorem{thm}{Théorème}
\newtheorem{prop}{Proposition}
\newtheorem{cor}{Corollaire}

\begin{document}
\title{Calculs effectifs de projections caractéristiques}
\author{J.~Puydt\\
Institut Joseph Fourier UMR5582\\
Grenoble, France\\
E-mail:\texttt{julien.puydt@ujf-grenoble.fr}
}

\maketitle

\begin{abstract}
  This article aims to make explicit the characteristic projection
  introduced by Panchishkin in \cite{Panchishkin2002b} and gives
  various results on how to calculate it effectively on a computer. It
  ends with concrete examples putting those ideas in application.
\end{abstract}

\section{Introduction}\label{sec:intr}

On souhaite dans cet article montrer comment rendre explicite et
effective informatiquement la méthode de la projection caractéristique
introduite par Pantchichkine dans \cite{Panchishkin2002b}, qui sert à
obtenir des congruences pour les valeurs spéciales des fonctions $L$
de formes modulaires. Pour cela, on va dans cette introduction
commencer par fixer quelques notations sur la théorie des formes
modulaires, et rappeler les énoncés principaux nécessaires. On
expliquera ensuite ce qu'est la projection caractéristique, puis on
énoncera les résultats proprement dits.

La suite de l'article consistera en des preuves et explications plus
précises et quelques exemples concrets de calculs (avec l'aide de
SAGE\footnote{version 5.4, disponible sur
  \texttt{http//www.sagemath.org}}).

Pour les définitions précises et les preuves, le lecteur est invité à
consulter des références comme les livres de Serre~\cite{Serre1970a},
de Miyake~\cite{Miyake1989a} ou de Diamond et
Shurman~\cite{DiamondShurman2005a}. On trouvera aussi un survol plus
détaillé et tourné vers les calculs effectifs dans le livre de
Stein~\cite{Stein2007a} : c'est en effet là que sont expliqués les
algorithmes utilisés par SAGE pour calculer les espaces de formes
modulaires et les opérateurs de Hecke associés.

Étant donnés $k\geqslant2$ et $N\geqslant1$, $\mathcal
M_k(\Gamma_1(N))$ dénote l'espace vectoriel des formes modulaires de
poids $k$ et de niveau $N$ pour $\Gamma_1(N)$, et $\mathcal
S_k(\Gamma_1(N))$ son sous-espace des formes paraboliques, sur lequel
on va se concentrer ; on peut voir un élément $f$ soit comme une
fonction $z\mapsto f(z)$ avec $z$ dans le demi-plan de Poincaré, soit
via le développement en série de Fourier à l'infini comme une fonction
$q\mapsto f(q)$ avec $q$ dans le disque unité, ce second point de vue
étant privilégié.

On sait en outre que ces espaces peuvent se décomposer en somme
directe de sous-espaces, indexée par les caractères de Dirichlet
$\chi$ de conducteur divisant $N$ :
\[
\left\{
\begin{array}{rcl}
  \mathcal M_k(\Gamma_1(N)) &=&\mathop{\oplus}_\chi\mathcal M_k(N,\chi) \\
  \mathcal S_k(\Gamma_1(N)) &=&\mathop{\oplus}_\chi\mathcal S_k(N,\chi) \\
\end{array}
\right.
\]

Entre ces espaces de formes modulaires, on sait définir des familles
d'opérateurs ; en particulier les opérateurs de Hecke,
$(T_n)_{n\geq2}$. La famille est engendrée par les opérateurs
$(T_p)_{p\ \mathrm{premier}}$, qui commutent deux à deux. Lorsque $p$
ne divise pas le niveau, $T_p$ est un endomorphisme de
l'espace. Lorsque $p$ divise le niveau, on parle d'opérateur de Hecke
spécial (et on écrit parfois $U_p$, ce que l'on fera ici) ; dans ce
cas, l'action sur les formes modulaires est un peu différente : elle
est décrite par exemple par le lemme~1 dans l'article de
Li~\cite{Li1975a}, et signifie intuitivement que l'opérateur fait
``disparaître'' une puissance de $p$ du niveau, sans pouvoir faire
disparaître complètement $p$. On verra plus loin une représentation
graphique de la situation : on obtient une tour d'espaces de formes
modulaires de dimensions finies, avec des inclusions croissantes avec
le niveau, mais une dimension qui croît assez rapidement ; la
projection est un moyen astucieux de contrer cette croissance.

Expliquons rapidement la construction telle qu'elle est préséntée
dans~\cite{Panchishkin2002b} : étant donnée une valeur propre $\alpha$
de $U_p$ agissant au rez-de-chaussée de la tour, on définit pour tout
$\mu\geqslant0$ le sous-espace caractéristique associé à
$(U_p,\alpha)$:
\[
\mathcal S_k^\alpha(Np^\mu,\chi)
=
\mathop{\cup}_{n\geqslant0}\mathrm{Ker}(U_p-\alpha\mathrm{Id})^n
\]
permettant de définir une projection $\pi_{\alpha,\mu}$ via le choix
d'un supplémentaire:
\[
\mathop{\cap}_{n\geqslant0}\mathrm{Im}(U_p-\alpha\mathrm{Id})^n
\]
la projection étant polynomiale en $U_p$ ; en faisant varier le niveau
par des puissances de $p$, on obtient alors un diagramme avec
projections $\pi_{\alpha,\_}$:
\[
\xymatrix{
\mathcal S_k(Np^\mu,\chi)
\ar[d]_{U_p}
\ar[r]^{\pi_{\alpha,\mu}}
&
\mathcal S_k^\alpha(Np^\mu,\chi)
\ar[d]^{U_p}
\\
\dots
\ar[d]_{U_p}&
\dots
\ar[d]^{U_p}
\\
\mathcal S_k(Np^2,\chi)
\ar[d]_{U_p}
\ar[r]^{\pi_{\alpha,2}}
&
\mathcal S_k^\alpha(Np^2,\chi)
\ar[d]^{U_p}
\\
\mathcal S_k(Np,\chi)
\ar@(dl,dr)[]_{U_p}
\ar[r]^{\pi_{\alpha,1}}
&
\mathcal S_k^\alpha(Np,\chi)
\ar@(dl,dr)[]_{U_p}
}
\]
dont la proposition 4.3 (de~\cite{Panchishkin2002b}, toujours) affirme
que les flèches verticales à droites sont des isomorphismes, d'où :
$\pi_{\alpha,\mu}=(U_p^\mu)^{-1}\pi_{\alpha,1}U_p^\mu$.

Ce dernier résultat est intéressant théoriquement, mais aussi
pratiquement : on est ramené à appliquer successivement $U_p$ sur des
espaces de plus en plus petits, puis projeter au plus bas niveau,
avant de finalement appliquer un inverse en petite dimension ; mieux,
d'après l'article \cite{ColemanEdixhoven1998a} de Coleman et
Edixhoven, on sait pour $k=2$ et conditionnellement à une conjecture
de Tate pour $k>2$ que l'opérateur $U_p$ est semi-simple au
rez-de-chaussée, donc l'inverse en question est une simple
multiplication par $\alpha^{-1}$. On en voit un exemple concret en
\ref{sec:calcul-tour-modulaire}.

Remarquons d'ailleurs que cela signifie aussi que dans les étages
supérieurs, tous les espaces caractéristiques associés à des valeurs
propres non nulles de $U_p$ sont réduits aux espaces propres : la
seule obstruction possible à la semi-simplicité de $U_p$ est la valeur
propre zéro. Donnons un exemple concret: pour $k=2, N=2, p=5$ et
$\chi$ modulo $10$ vérifiant $\chi(7)=1$ ou $\chi(7)=-1$, l'espace est
réduit à zéro en niveau $Np$ donc il n'y a pas de valeur propre non
nulle, est un plan en niveau $Np^2$ sur lequel $U_p$ agit par zéro, et
est de dimension $28$ en niveau $Np^3$ ; $U_p$ est alors une matrice
$28\times28$ de rang $2$, dont une forme réduite de Jordan comporte
deux blocs $\begin{pmatrix}0&1\\0&0\end{pmatrix}$ et vingt-quatre
blocs $\begin{pmatrix}0\end{pmatrix}$ -- ce qui montre que les
conditions de semi-simplicité de \cite{ColemanEdixhoven1998a} sont
optimales.

On prouve ici tout d'abord le résultat élémentaire suivant:
\begin{prop}\label{prop:explicite}
  La projection caractéristique est un polynôme construit
  explicitement à partir du polynôme minimal de l'opérateur.
\end{prop}
(l'algorithme, utilisant l'algorithme d'Euclide étendu, est présenté
de façon théorique en section~\ref{sec:defin-elem}, puis en détaillant
plus les étapes algorithmiques en section~\ref{sec:calcul-direct})

Les deux résultats d'indépendance suivants le rendent plus utile :
\begin{prop}
\begin{itemize}
\item\label{prop:plusieurs_poly} Il existe toute une famille de
  polynômes convenables, qui donnent tous le même opérateur de
  projection (voir \ref{sec:defin-elem}).
\item\label{prop:indep_poly} L'opérateur de projection ne dépend pas
  du choix du polynôme annulateur choisi (voir
  \ref{sec:chang-de-polyn}).
\end{itemize}
\end{prop}

Le premier point peut permettre de simplifier les calculs, au prix
d'un polynôme de degré plus élevé (voir par exemple en
section~\ref{sec:calcul-sans-division} une variante de l'algorithme
avec moins de divisions) ; le second permet d'initier le calcul avec
le polynôme minimal ou le polynôme caractéristique par exemple, ce qui
donne un peu plus de liberté dans les calculs -- mais ne les simplifie
que de façon marginale ; le résultat le plus important pour les
applications, qui permet des gains conséquents, est celui-ci:
\begin{thm}\label{thm:equivariance}
  La contruction de la projection canonique est équivariante ; plus
  précisément, étant donnés deux espaces vectoriels de dimension finie
  $E$ et $F$, $u\in\mathcal L(E)$ et $v\in\mathcal L(F)$,
  $j\in\mathcal L(E,F)$ et $\alpha\in\CC$, tels que $j\circ u=v\circ
  j$, alors $j\circ\pi_{E,\alpha}(u)=\pi_{F,\alpha}(v)\circ j$.
\end{thm}
(il est prouvé en \ref{sec:stab-par-morph})

En effet, ce résultat a comme corollaire immédiat :
\begin{cor}\label{cor:equiv_chars}
  Dans l'égalité suivante, où $p$ divise $N$, étant donné
  $\alpha\in\CC$ :
  \[\mathcal S_k(\Gamma_1(N)) =\mathop{\oplus}_\chi\mathcal S_k(N,\chi) \]
  on peut calculer les opérateurs de projection associés à
  $(U_p,\alpha)$ sur les sous-espaces de droite (plus petits), et
  obtenir par recollement l'opérateur de projection sur l'espace de
  gauche.
\end{cor}

Pour énoncer le second corollaire ci-dessous, il faut en dire plus sur
les espaces de formes modulaires et expliquer un peu la théorie de
Atkin-Lehner-Li : en plus des opérateurs de Hecke, on peut considérer
les opérateurs $V_d$, dont la définition est assez simple: pour
$d\geqslant1$, on a $(V_df)(z)=f(d~z)$ ou $(V_df)(q)=f(q^d)$, ce qui
revient au même. En particulier, $V_1$ est l'inclusion naturelle
mentionnée précédemment.

Si on fixe un niveau $N\geqslant1$, on appelle forme ancienne de
$\mathcal S_k(\Gamma_1(N))$ toute combinaison linéaire de formes qui
proviennent d'espaces de niveaux strictement inférieur via les
opérateurs $V$ -- cela définit le sous espace $\mathcal
S_k(\Gamma_1(N))_a$. Le supplémentaire orthogonal de ce sous-espace
(vis-à-vis du produit scalaire de Petersson) est appelé espace des
nouvelles formes, noté $\mathcal S_k(\Gamma_1(N))_n$. La théorie
d'Atkin-Lehner-Li (développée par Atkin et Lehner dans
\cite{AtkinLehner1970a}, puis par Li dans \cite{Li1975a}) a pour
conséquence que :
\[
\mathcal S_k(\Gamma_1(N)=\oplus_{dM|N}V_d\mathcal S_k(\Gamma_1(M))_n
\]
et cette décomposition est par exemple utilisée pour les calculs
d'espaces modulaires par informatique dans SAGE, comme expliqué
dans le livre de Stein~\cite{Stein2007a}.

On peut maintenant énoncer le second corollaire du théorème
d'équivariance:
\begin{cor}\label{cor:equiv_ALL}
  Dans la décomposition précédente, si les seuls $d$ qui interviennent
  réellement (pour lesquels les sous-espaces ne sont pas nuls) sont
  premiers avec $p$, alors on peut calculer les opérateurs de
  projections sur les sous-espaces de droite (de dimensions
  inférieures) et obtenir l'opérateur de projection sur l'espace de
  gauche par recollement.
\end{cor}
(l'exemple \ref{sec:calcul-ALL} illustre ce résultat)

Je souhaite remercier Alexei Pantchichkine pour son soutien
indéfectible, et William Stein pour m'avoir donné accès aux machines
du réseau math.washington.edu\footnote{Les calculs présentés dans cet
  article sont tous aisés et rapides sur une machine simple, mais les
  tâtonnements et expériences nécessaires à leur recherche ont parfois
  nécessité l'accès à une puissance de calcul plus conséquente ; les
  machines en question ont été financées par "National Science
  Foundation Grant No. DMS-0821725".}.

\section{Généralités d'algèbre linéaire}\label{sec:gener-dalg-line}

\subsection{Définition élémentaire}\label{sec:defin-elem}

On se donne $E$ un $\KK$-espace vectoriel de dimension finie et $u$ un
endomorphisme de $E$, admettant une valeur propre $\alpha$.

Le polynôme minimal de $u$ sur $E$ admet une factorisation de la forme
$(X-\alpha)^\nu Q$ où $\nu\geqslant1$ et $Q(\alpha)\neq0$ ; d'après le
théorème de Bézout, il existe alors un couple de polynômes $(A,B)$ tel
que $(X-\alpha)^\nu A+QB=1$, dont on déduit une décomposition de $E$
en supplémentaires: $E=\Ker(u-\alpha)^\nu\oplus\Ker~Q(u)$.

L'espace $\Ker(u-\alpha)^\nu$ est appelé espace caractéristique de $u$
associé à la valeur propre $\alpha$ ; c'est le plus grand espace sur
lequel $u-\alpha$ est nilpotent. La décomposition en supplémentaires
stables par $u$ précédente permet alors de définir une projection sur
cet espace caractéristique ; on vérifie aisément que $Q(u)B(u)$ donne
une expression de l'opérateur associé. Cela établit déjà la
proposition~\ref{prop:explicite}.

Une remarque simple mais fondamentale est que le couple de polynômes
$(A,B)$ considéré n'est pas unique, mais que deux couples admissibles
distincts diffèrent d'un multiple de $(X-\alpha)^\nu Q$ -- le polynôme
annulateur de $u$. Donc l'opérateur $Q(u)B(u)$ ne dépend pas du couple
choisi ; ce qui prouve la première
partie de la proposition~\ref{prop:plusieurs_poly}.

On notera $\pi_{E,\alpha}(u)$ l'opérateur de projection
caractéristique de $u$ associé à $\alpha$ sur l'espace $E$.

\subsection{Stabilité par extension}\label{sec:stab-par-extens}

On se donne maintenant un $\KK$-espace vectoriel $E$ de dimension
finie, $u$ un endomorphisme de $E$ ; et on suppose que $F$ est un
sous-espace vectoriel de $E$ stable par $u$, et que $\alpha$ est une
valeur propre de $u$ sur $F$.

On dispose donc de deux opérateurs de projection caractéristique,
$\pi_{E,\alpha}(u)$ et $\pi_{F,\alpha}(u)$, que l'on souhaite comparer
sur $F$.

Les polynômes minimaux de $u$ sur $F$ et $E$ sont de la forme
respectivement $(X-\alpha)^\nu Q$ et $(X-\alpha)^\mu R$, avec
$1\leqslant\mu\leqslant\nu$, $Q$ diviseur de $R$ et
$R(\alpha)\neq0$. Écrivons $R=QS$ pour fixer les notations.

L'application du théorème de Bézout pour le polynôme minimal de $u$
sur $E$ fournit un couple $(C,D)$ de polynômes tel que $(X-\alpha)^\nu
C+RD=1$ ; on a donc $\pi_{E,\alpha}(u)=R(u)D(u)$.

Mais cette même expression peut s'écrire
$(X-\alpha)^\mu(X-\alpha)^{\nu-\mu}C+QSD=1$ ; le couple de polynômes
$((X-\alpha)^{\nu-\mu}C, SD)$ est donc admissible pour l'application
du théorème de Bézout pour le polynôme minimal de $u$ sur $F$ ; on a
donc $\pi_{F,\alpha}(u)=Q(u)(SD)(u)$.

Il reste à écrire $Q(u)(SD)(u)=(QS)(u)D(u)=R(u)D(u)$ pour conclure que
l'on a bien $\pi_{F,\alpha}(u)=\pi_{E,\alpha}(u)_{|F}$.

Ce résultat est un lemme de base pour la preuve de la
proposition~\ref{prop:indep_poly} et du
théorème~\ref{thm:equivariance}.

\subsection{Changement de polynôme annulateur}\label{sec:chang-de-polyn}

Une première variante importante de ce raisonnement est le cas où l'on
ne considère $u$ que sur un espace, mais avec deux polynômes
annulateurs. Comme le polynôme minimal est un diviseur commun à deux
tels polynômes, il suffit de considérer la situation d'un polynôme
annulateur quelconque par rapport au polynôme minimal : une
compatibilité dans ce cas montre alors que la compatibilité est
générale.

Maintenant, si compare avec la preuve précédente, on voit que la seule
chose que l'on ait utilisée sur le polynôme minimal sur le sur-espace
est qu'il était un multiple du polynôme minimal sur l'espace de base ;
la même preuve s'applique donc : le calcul de la projection
caractéristique ne dépend donc pas du polynôme annulateur considéré ;
cela établit la seconde partie de la
proposition~\ref{prop:indep_poly}.

\subsection{Stabilité par morphisme}\label{sec:stab-par-morph}

Une seconde variante importante de ce raisonnement est celui où l'on
considère deux espaces vectoriels de dimensions finies $E$ et $F$, $u$
un endomorphisme de $E$, $v$ un endomorphisme de $F$ et $j$ une
application linéaire injective de $E$ dans $F$ telle que $j\circ
u=v\circ j$. Si $\alpha$ est une valeur propre de $u$ sur $E$, c'est
aussi une valeur propre de $v$ sur $F$ et on souhaite comparer
$\pi_{E,\alpha}(u)$ et $\pi_{F,\alpha}(v)$ ; la situation peut se
visualiser ainsi :
\[
\xymatrix{
E
\ar@(ul,dl)[]_u
\ar[r]^j
&
F
\ar@(ur,dr)[]^v
}
\]

Si $P$ est un polynôme annulateur de $u$ sur $E$, alors la relation
$jP(u)=P(v)j$ montre que $P$ est un polynôme annulateur de $v$ sur
$j(E)$ ; les mêmes raisonnements que précédemment montrent que le
calcul de la projection caractéristique peut se faire au départ comme
à l'arrivée :
\[
j\circ\pi_{E,\alpha}(u)=\pi_{F,\alpha}(v)\circ j
\]

C'est le théorème d'équivariance~\ref{thm:equivariance}.

\section{Pratique}

\subsection{Calcul direct}\label{sec:calcul-direct}

Il est assez facile de déterminer le polynôme minimal d'un
endomorphisme donné sous forme matricielle via l'algorithme de
Wiedemann, que l'on trouve détaillé en section~7.5.3 du livre de
Stein~\cite{Stein2007a}.

Il est aisé de calculer un couple explicite dans le théorème de Bézout
via l'algorithme d'Euclide étendu, tel qu'il est présenté par exemple
en section~3.2 du livre de Cohen~\cite{Cohen1993a}.

Le calcul d'un polynôme d'endomorphisme à partir d'un polynôme
explicite et d'un endomorphisme donné par exemple sous forme d'une
matrice, est aussi facile, via l'algorithme de Horner par exemple.

Finalement, en combinant ces étapes, on obtient facilement une
expression matricielle de la projection canonique ; la principale gêne
dans les calculs est le passage au corps de nombres $\QQ(\alpha)$,
dans lequel on réalise le calcul de divisions euclidiennes pour
l'algorithme d'Euclide.

Un cas particulier important est celui où l'on sait que $\alpha$ est
une racine simple du polynôme annulateur considéré ; en effet,
l'algorithme d'Euclide s'applique en une seule division euclidienne,
que l'on peut calculer par la variante de l'algorithme de Horner qui
traite la division par $X-\alpha$ ; elle ne fait donc intervenir que
peu de produits et de sommes et une seule division.

Un sous-cas intéressant du point précédent est celui de la tour
modulaire au rez-de-chaussée de laquelle on sait que l'opérateur $U_p$
est semi-simple (sous couvert d'une conjecture de Tate pour $k>2$),
car alors on peut simplifier le polynôme annulateur choisi pour se
ramener à une racine simple.

\subsection{Calcul sans division intermédiaire}\label{sec:calcul-sans-division}

On sait que les divisions du calcul précédent sont incontournables en
général.

Cependant, les polynômes avec lesquels on travaille ont une forme
assez particulière ; en utilisant cette spécificité, on peut
contourner en partie cette difficulté de la façon suivante : étant
donnée une factorisation d'un polynôme annulateur sous la forme
$(X-\alpha)^\nu P$, avec $P(\alpha)=a\neq0$, on écrit
$P=a+(X-\alpha)^\mu Q$ avec $\mu\geqslant1$.

Pour $e\geqslant1$, on a alors:
\[
P
\prod_{j=1}^e\left(a^{2^{j-1}}+(-1)^j(X-\alpha)^{2^{j-1}\mu}Q^{2^{j-1}}\right)
=
a^{2^e}-(X-\alpha)^{2^e\mu}Q^{2^e}
\]
de cette façon, si $2^e\mu\geqslant\nu$, ce qui se produit pour $e$
assez petit, on a :
\[
P
\prod_{j=1}^e\left(a^{2^{j-1}}+(-1)^j(X-\alpha)^{2^{j-1}\mu}Q^{2^{j-1}}\right)
+(X-\alpha)^\nu(X-\alpha)^{2^e\mu-\nu}Q^{2^e}
=
a^{2^e}
\]

expression qu'il suffit de diviser une fois par le scalaire $a^{2^e}$
pour écrire un couple admissible. On remplace le coût des divisions
euclidiennes par un calcul avec un polynôme de degré plus élevé ;
c'est un compromis.

Remarquons aussi qu'une bibliothèque de calcul en théorie des nombres
comme FLINT\footnote{\texttt{http://www.flintlib.org/}} dispose de
routines optimisées pour les opérations de décalage (\og Taylor
shift\fg{}) et de mise au carré, qui peuvent être mises à profit pour
les calculs précédents.

\subsection{Détermination des dimensions}

On a vu précédemment comment calculer très explicitement et
complètement les opérateurs de projections caractéristiques, et on
sait qu'alors la dimension de l'espace sur lequel on projette
s'obtient par un simple calcul de trace.

On peut cependant souhaiter connaître la dimension de l'espace sans
avoir calculé la projection, par exemple parce que l'on n'a pas encore
choisi quelle valeur propre utiliser. Dans ce cas, il suffit de
calculer le polynôme caractéristique de l'opérateur : en effet, il
suffira alors de le factoriser pour obtenir les ordres de multiplicité
algébrique, donc les dimensions des espaces caractéristiques associés.

\section{Exemples concrets}

\subsection{Calcul via la théorie d'Atkin-Lehner-Li}\label{sec:calcul-ALL}

On vérifie aisément (avec SAGE, qui permet de déterminer les espaces
réduits à $\{0\}$) que l'on a la décomposition suivante:
\[
\mathcal S_2(\Gamma_1(30))
=
V_1\mathcal S_2(\Gamma_1(15))_{\mathrm{new}}
\oplus
V_2\mathcal S_2(\Gamma_1(15))_{\mathrm{new}}
\oplus
V_1\mathcal S_2(\Gamma_1(30))_{\mathrm{new}}
\]
(décomposition dans laquelle les espaces sont de dimensions respectives 9, 1, 1 et 7).

Une base adaptée à cette décomposition est :
\[
\begin{array}{c}
q - q^2 - q^3 - q^4 + q^5 + q^6 + 3q^8 + q^9 +O( q^{10})\\
q^2 - q^4 - q^6 - q^8 + O(q^{10}),\\
q + 3q^8 - 14q^9 +O(q^{10}),\\
q^2 + 2q^8 - 10q^9 +O(q^{10}),\\
q^3 + 2q^8 - 9q^9 +O(q^{10}),\\
q^4 + q^8 - 5q^9 +O(q^{10}),\\
q^5 + q^8 - 5q^9 +O(q^{10}),\\
q^6 - 2q^9 +O(q^{10}),\\
q^7 + q^8 - 3q^9 +O(q^{10})\\
\end{array}
\]

La base calculée par SAGE directement pour le même espace est une base
de Miller ; donc très simple:
\[
\begin{array}{ll}
q+O(q^{10}),
q^2+O(q^{10}),
q^3+O(q^{10}),
q^4+O(q^{10}),
q^5+O(q^{10}),\\
q^6+O(q^{10}),
q^7+O(q^{10}),
q^8+O(q^{10}),
q^9+O(q^{10})
\end{array}
\]

La matrice de passage est donc:
\[
P
=
\left(\begin{array}{rrrrrrrrr}
1 & -1 & -1 & -1 & 1 & 1 & 0 & 3 & 1 \\
0 & 1 & 0 & -1 & 0 & -1 & 0 & -1 & 0 \\
1 & 0 & 0 & 0 & 0 & 0 & 0 & 3 & -14 \\
0 & 1 & 0 & 0 & 0 & 0 & 0 & 2 & -10 \\
0 & 0 & 1 & 0 & 0 & 0 & 0 & 2 & -9 \\
0 & 0 & 0 & 1 & 0 & 0 & 0 & 1 & -5 \\
0 & 0 & 0 & 0 & 1 & 0 & 0 & 1 & -5 \\
0 & 0 & 0 & 0 & 0 & 1 & 0 & 0 & -2 \\
0 & 0 & 0 & 0 & 0 & 0 & 1 & 1 & -3
\end{array}\right)
\]

L'endomorphisme $T_3$ agit par $-1$ sur les deux premiers sous-espaces, et par:
\[
\left(\begin{array}{rrrrrrr}
0 & 0 & -14 & 2 & 19 & 7 & -3 \\
0 & 0 & -10 & 1 & 15 & 4 & -3 \\
1 & 0 & -9 & 1 & 11 & 3 & -3 \\
0 & 0 & -5 & 0 & 8 & 3 & -3 \\
0 & 0 & -5 & 1 & 6 & 2 & 0 \\
0 & 1 & -2 & 0 & 2 & 0 & -2 \\
0 & 0 & -3 & 0 & 4 & 2 & 0
\end{array}\right)
\]
sur le dernier ; son polynôme caractéristique est alors:
$X^7+3X^6+5X^5+7X^4+X^3-5X^2-3X-9$. Il a une racine simple réelle,
$1$, deux racines simples conjuguées dans $\QQ(i)$ ($i$ et $-i$) et
quatre racines simples conjuguées dans un corps de nombres de degré
$4$ ; il s'agit donc aussi du polynôme minimal de l'endomorphisme.

Choisissons de calculer la projection caractéristique associée à la
valeur propre $1$. Sur les deux droites, elle est donnée par zéro ; et
sur le troisième espace, par le polynôme
$\frac1{68}(X^6+4X^5+9X^4+16X^3+17X^2+12X+9)$ appliqué à la matrice
précédente, soit :
\[
\left(\begin{array}{rrrrrrr}
-\frac{1}{2} & \frac{1}{2} & -\frac{1}{2} & -\frac{1}{2} & \frac{1}{2} & \frac{1}{2} & 2 \\
-\frac{1}{2} & \frac{1}{2} & -\frac{1}{2} & -\frac{1}{2} & \frac{1}{2} & \frac{1}{2} & 2 \\
-\frac{1}{4} & \frac{1}{4} & -\frac{1}{4} & -\frac{1}{4} & \frac{1}{4} & \frac{1}{4} & 1 \\
0 & 0 & 0 & 0 & 0 & 0 & 0 \\
-\frac{1}{4} & \frac{1}{4} & -\frac{1}{4} & -\frac{1}{4} & \frac{1}{4} & \frac{1}{4} & 1 \\
0 & 0 & 0 & 0 & 0 & 0 & 0 \\
-\frac{1}{4} & \frac{1}{4} & -\frac{1}{4} & -\frac{1}{4} & \frac{1}{4} & \frac{1}{4} & 1
\end{array}\right)
\]

Les résultats du paragraphe précédent ont été obtenus de la façon
suivante dans SAGE:
\begin{verbatim}
# Détermination de l'espace et de l'opérateur
S30=ModularForms(Gamma1(30),2).cuspidal_subspace()
S30n=S30.new_subspace()
T=S30n.T(3)
# Calcul du polynôme caractéristique
R.<t>=PolynomialRing(QQ)
P=T.charpoly(t)
# Calcul du polynôme à utiliser pour le calcul de la projection
Q=P.quo_rem(t-1)[0]
Pproj=xgcd(Q,t-1)[1]*Q
# Obtention de la matrice de la projection
Pproj(T.matrix())
\end{verbatim}

La matrice de la projection canonique sur l'espace complet est alors,
dans la base adaptée à la décomposition donnée ci-dessus:
\[
\left(\begin{array}{r|r|rrrrrrr}
0 & 0 & 0 & 0 & 0 & 0 & 0 & 0 & 0 \\
\hline
 0 & 0 & 0 & 0 & 0 & 0 & 0 & 0 & 0 \\
\hline
 0 & 0 & -\frac{1}{2} & \frac{1}{2} & -\frac{1}{2} & -\frac{1}{2} & \frac{1}{2} & \frac{1}{2} & 2 \\
0 & 0 & -\frac{1}{2} & \frac{1}{2} & -\frac{1}{2} & -\frac{1}{2} & \frac{1}{2} & \frac{1}{2} & 2 \\
0 & 0 & -\frac{1}{4} & \frac{1}{4} & -\frac{1}{4} & -\frac{1}{4} & \frac{1}{4} & \frac{1}{4} & 1 \\
0 & 0 & 0 & 0 & 0 & 0 & 0 & 0 & 0 \\
0 & 0 & -\frac{1}{4} & \frac{1}{4} & -\frac{1}{4} & -\frac{1}{4} & \frac{1}{4} & \frac{1}{4} & 1 \\
0 & 0 & 0 & 0 & 0 & 0 & 0 & 0 & 0 \\
0 & 0 & -\frac{1}{4} & \frac{1}{4} & -\frac{1}{4} & -\frac{1}{4} & \frac{1}{4} & \frac{1}{4} & 1
\end{array}\right)
\]

En calculant directement avec l'opérateur sur l'espace total, il faut
appliquer sa matrice au polynôme
$\frac1{272}X^8+6X^7+18X^6+38X^5+58X^4+62X^3+50X^2+30X+9$, d'où:
\[
\left(\begin{array}{rrrrrrrrr}
\frac{21}{8} & -\frac{21}{8} & \frac{21}{8} & \frac{21}{8} & -\frac{21}{8} & -\frac{21}{8} & -\frac{21}{2} & -\frac{21}{8} & \frac{21}{8} \\
\frac{7}{4} & -\frac{7}{4} & \frac{7}{4} & \frac{7}{4} & -\frac{7}{4} & -\frac{7}{4} & -7 & -\frac{7}{4} & \frac{7}{4} \\
\frac{7}{4} & -\frac{7}{4} & \frac{7}{4} & \frac{7}{4} & -\frac{7}{4} & -\frac{7}{4} & -7 & -\frac{7}{4} & \frac{7}{4} \\
\frac{9}{8} & -\frac{9}{8} & \frac{9}{8} & \frac{9}{8} & -\frac{9}{8} & -\frac{9}{8} & -\frac{9}{2} & -\frac{9}{8} & \frac{9}{8} \\
\frac{7}{8} & -\frac{7}{8} & \frac{7}{8} & \frac{7}{8} & -\frac{7}{8} & -\frac{7}{8} & -\frac{7}{2} & -\frac{7}{8} & \frac{7}{8} \\
\frac{1}{2} & -\frac{1}{2} & \frac{1}{2} & \frac{1}{2} & -\frac{1}{2} & -\frac{1}{2} & -2 & -\frac{1}{2} & \frac{1}{2} \\
\frac{3}{8} & -\frac{3}{8} & \frac{3}{8} & \frac{3}{8} & -\frac{3}{8} & -\frac{3}{8} & -\frac{3}{2} & -\frac{3}{8} & \frac{3}{8} \\
\frac{1}{8} & -\frac{1}{8} & \frac{1}{8} & \frac{1}{8} & -\frac{1}{8} & -\frac{1}{8} & -\frac{1}{2} & -\frac{1}{8} & \frac{1}{8} \\
\frac{1}{4} & -\frac{1}{4} & \frac{1}{4} & \frac{1}{4} & -\frac{1}{4} & -\frac{1}{4} & -1 & -\frac{1}{4} & \frac{1}{4}
\end{array}\right)
\]

matrice qui est équivalente à celle trouvée précédemment, à $P$ près.

\subsection{Calcul avec une tour modulaire}\label{sec:calcul-tour-modulaire}

Cet exemple est un cas où les calculs sont encore possibles à tous les
niveaux, mais illustre le gain du passage à une dimension inférieure.

On choisit un niveau de base $30$ et le nombre premier $3$, et un
caractère de Dirichlet trivial, de sorte que la situation est:
\[
\xymatrix{
(\mathrm{dimension}\ 11) &
\mathcal S_2(90,1)
\ar[r]^{\pi_{-1,1}}
\ar[d]_{U_3}
&
\mathcal S_2^{-1}(90,1)
\ar[d]^{U_3}
\\
(\mathrm{dimension}\ 3) &
\mathcal S_2(30,1)
\ar[r]_{\pi_{-1,0}}
\ar@(dl,dr)[]_{U_3}
&
\mathcal S_2^{-1}(30,1)
\ar@(dl,dr)[]_{U_3}
}
\]

On constate simplement que $U_3$ agit par simple multiplication par
$-1$ sur les espaces caractéristiques, qui sont des plans ; la tour
permet donc de calculer la projection en haut via un calcul en bas,
en dimension inférieure.

Ici, la projection au rez-de-chaussée se calcule avec le polynôme
$-\frac{1}{4} X^{2} - \frac{1}{2} X + \frac{3}{4}$ et est donnée par
la matrice:
\[
\left(\begin{array}{rrr}
\frac{1}{2} & \frac{1}{2} & -\frac{1}{2} \\
0 & 1 & 0 \\
-\frac{1}{2} & \frac{1}{2} & \frac{1}{2}
\end{array}\right)
\]

À l'étage au dessus, elle est donnée par le polynôme $-\frac{17}{4}
X^{10} - \frac{1}{2} X^{9} + \frac{19}{4} X^{8}$, et la matrice:
\[
\left(\begin{array}{rrrrrrrrrrr}
0 & 0 & 0 & 0 & 0 & 0 & 0 & 0 & 0 & 0 & 0 \\
0 & 0 & 0 & 0 & 0 & 0 & 0 & 0 & 0 & 0 & 0 \\
-\frac{1}{2} & -\frac{1}{2} & \frac{1}{2} & \frac{3}{2} & -\frac{1}{2} & \frac{1}{2} & 0 & -\frac{1}{2} & -\frac{1}{2} & -\frac{1}{2} & 2 \\
0 & 0 & 0 & 0 & 0 & 0 & 0 & 0 & 0 & 0 & 0 \\
0 & 0 & 0 & 0 & 0 & 0 & 0 & 0 & 0 & 0 & 0 \\
0 & -1 & 0 & 1 & 0 & 1 & 0 & 1 & 0 & -1 & 0 \\
0 & 0 & 0 & 0 & 0 & 0 & 0 & 0 & 0 & 0 & 0 \\
0 & 0 & 0 & 0 & 0 & 0 & 0 & 0 & 0 & 0 & 0 \\
\frac{1}{2} & -\frac{1}{2} & -\frac{1}{2} & -\frac{1}{2} & \frac{1}{2} & \frac{1}{2} & 0 & \frac{3}{2} & \frac{1}{2} & -\frac{1}{2} & -2 \\
0 & 0 & 0 & 0 & 0 & 0 & 0 & 0 & 0 & 0 & 0 \\
0 & 0 & 0 & 0 & 0 & 0 & 0 & 0 & 0 & 0 & 0
\end{array}\right)
\]

En choisissant une forme modulaire $f$ dans l'espace $\mathcal
S_2(90,1)$, on vérifie aisément l'égalité
$\pi_{-1,1}(f)=-\pi_{-1,0}U_3(f)$ ; la session ci-après utilise les
deux fonctions accessoires données en section~\ref{sec:code-sage} :
\begin{verbatim}
sage: # Chargement des fonctions d'aide
sage: load('/home/jpuydt/Recherche/abaisse.sage')
sage: load('/home/jpuydt/Recherche/polyproj.sage')
sage: # On travaille avec des polynômes, pas des expressions
sage: R.<X>=PolynomialRing(QQ)
sage: # Données du calcul
sage: N=10
sage: p=3
sage: N0=N*p
sage: # Choix du caractère
sage: chi=DirichletGroup(N0)[0]
sage: # Calcul de la projection en bas
sage: Eb=ModularForms(chi).cuspidal_subspace()
sage: Ub=Eb.T(p)
sage: polb=polyproj(Ub.charpoly(X),-1,2)
sage: # Calcul de la projection en haut
sage: Eh=ModularForms(chi.extend(N0*p)).cuspidal_subspace()
sage: Uh=Eh.T(p)
sage: polh=polyproj(Uh.charpoly(X),-1,2)
sage: # Comparaison des deux
sage: f=Eh.gens()[2]
sage: g=abaisse(Eb, Uh(f))
sage: ordre=Eh.sturm_bound ()
sage: polb(Ub)(g).qexp(ordre) + polh(Uh)(f).qexp(ordre)
O(q^37)
sage: # On obtient bien zéro!
\end{verbatim}

\section{Code SAGE}\label{sec:code-sage}

\subsection{Fonction \texttt{polyproj}}

\begin{verbatim}
def polyproj(poly, racine, ordre):
    """Reçoit en argument un polynôme annulateur d'un endomorphisme,
    une de ses racines et son ordre et retourne un polynôme permettant
    de calculer la projection caractéristique sur la racine
    correspondante"""
    var=poly.parent().gen()
    pow=(var-racine)**ordre 
    Q=poly.quo_rem(pow)[0]
    return Q*xgcd(Q,pow)[1]
\end{verbatim}

\subsection{Fonction \texttt{abaisse}}

\begin{verbatim}
def abaisse (espace, forme):
    """Reçoit en argument un espace de formes modulaires et une forme
    modulaire a priori dans un espace de niveau plus élevé, mais dont
    on sait qu'elle devrait être plus bas, et retourne la 'même' forme
    en bas."""
    coeffs=espace.find_in_space(forme)
    return espace.linear_combination_of_basis(coeffs)
\end{verbatim}

\bibliography{../Abords/bibliography}

\begin{thebibliography}{1}
\expandafter\ifx\csname fonteauteurs\endcsname\relax
\def\fonteauteurs{\scshape}\fi

\bibitem{AtkinLehner1970a}
A.~O.~L. \bgroup\fonteauteurs\bgroup Atkin\egroup\egroup{} et
  J.~\bgroup\fonteauteurs\bgroup Lehner\egroup\egroup{} :
\newblock {H}ecke operators on {$\Gamma_0(m)$}.
\newblock {\em {M}athematische {A}nnalen}, 185\string:\penalty500\relax
  134--160, 1970.

\bibitem{Cohen1993a}
H.~\bgroup\fonteauteurs\bgroup Cohen\egroup\egroup{} :
\newblock {\em {A} course in computational algebraic number theory}, volume 138
  de {\em {G}raduate texts in mathematics}.
\newblock {S}pringer-{V}erlag, 1993.

\bibitem{ColemanEdixhoven1998a}
R.~F. \bgroup\fonteauteurs\bgroup Coleman\egroup\egroup{} et S.~J.
  \bgroup\fonteauteurs\bgroup Edixhoven\egroup\egroup{} :
\newblock {O}n the semi-simplicity of the {$U_p$}-operator on modular forms.
\newblock {\em {M}athematische {A}nnalen}, 310\string:\penalty500\relax
  119--127, 1998.

\bibitem{DiamondShurman2005a}
F.~\bgroup\fonteauteurs\bgroup Diamond\egroup\egroup{} et
  J.~\bgroup\fonteauteurs\bgroup Shurman\egroup\egroup{} :
\newblock {\em {A} first course in modular forms}, volume 228 de {\em
  {G}raduate texts in mathematics}.
\newblock {S}pringer-{V}erlag, 2005.

\bibitem{Li1975a}
W.-C.~W. \bgroup\fonteauteurs\bgroup Li\egroup\egroup{} :
\newblock {N}ewforms and functional equations.
\newblock {\em {M}athematische {A}nnalen}, 212(4)\string:\penalty500\relax
  285--315, 1975.

\bibitem{Miyake1989a}
T.~\bgroup\fonteauteurs\bgroup Miyake\egroup\egroup{} :
\newblock {\em {M}odular forms}.
\newblock {S}pringer, 1989.

\bibitem{Panchishkin2002b}
A.~A. \bgroup\fonteauteurs\bgroup Panchishkin\egroup\egroup{} :
\newblock {A} new method of constructing {$p$}-adic {$L$}-functions associated
  with modular forms.
\newblock {\em {M}oscow mathematical journal}, 2(2)\string:\penalty500\relax
  313--328, 2002.

\bibitem{Serre1970a}
J.-P. \bgroup\fonteauteurs\bgroup Serre\egroup\egroup{} :
\newblock {\em {C}ours d'arithmétique}.
\newblock {P}resses universitaires de {F}rance, 1970.

\bibitem{Stein2007a}
W.~A. \bgroup\fonteauteurs\bgroup Stein\egroup\egroup{} :
\newblock {\em {M}odular forms, a computational approach}, volume~79 de {\em
  {G}raduate studies in mathematics}.
\newblock {A}merican mathematical society, 2007.

\end{thebibliography}
\bibliographystyle{plain-fr}

\end{document}